\documentclass[12pt,a4paper]{amsart}

\usepackage[colorlinks,pagebackref,hypertexnames=false]{hyperref}

\usepackage[a4paper,left=1.8cm,right=1.8cm,top=2.5cm,bottom=2.5cm]{geometry}

\makeatletter
\renewcommand\normalsize{%
    \@setfontsize\normalsize{11.7}{14pt plus .3pt minus .3pt}%
    \abovedisplayskip 10\p@ \@plus4\p@ \@minus4\p@
    \abovedisplayshortskip 6\p@ \@plus2\p@
    \belowdisplayshortskip 6\p@ \@plus2\p@
    \belowdisplayskip \abovedisplayskip}
\renewcommand\small{%
    \@setfontsize\small{9.5}{12\p@ plus .2\p@ minus .2\p@}%
    \abovedisplayskip 8.5\p@ \@plus4\p@ \@minus1\p@
    \belowdisplayskip \abovedisplayskip
    \abovedisplayshortskip \abovedisplayskip
    \belowdisplayshortskip \abovedisplayskip}
\renewcommand\footnotesize{%
    \@setfontsize\footnotesize{8.5}{9.25\p@ plus .1pt minus .1pt}
    \abovedisplayskip 6\p@ \@plus4\p@ \@minus1\p@
    \belowdisplayskip \abovedisplayskip
    \abovedisplayshortskip \abovedisplayskip
    \belowdisplayshortskip \abovedisplayskip}
\ifdefined\pdfpagewidth
\else
\fi
\calclayout
\makeatother

\usepackage{amsmath, amsthm, amsfonts, amssymb, mathtools, commath, upgreek}
\usepackage{graphicx}
\usepackage[english]{babel}
\usepackage{hyperref}
\usepackage{microtype}
\usepackage{enumitem}

\usepackage{lmodern}
\usepackage{bbm}
\usepackage{float}
\usepackage{xurl}
\usepackage{ragged2e}

\usepackage{empheq, mdframed}
\usepackage[all]{xy}

\usepackage{etoolbox}
\AtBeginEnvironment{thebibliography}{%
  \setlength{\itemsep}{0pt}%
  \setlength{\parskip}{0pt}%
  \setlength{\topsep}{0pt}%
  \setlength{\partopsep}{0pt}%
}

\newlist{defenum}{enumerate}{3}
\setlist[defenum]{label=(\alph*),leftmargin=*,align=left}
\newlist{defsubenum}{enumerate}{1}
\setlist[defsubenum]{label=(\roman*),leftmargin=*,align=left}

\usepackage{enumitem}
\setlist{nosep}
\setlist[enumerate]{label=(\roman*)}

\hypersetup{
    colorlinks=true,
    linkcolor=red,
    citecolor=blue,
    filecolor=blue,
    urlcolor=blue
}

\usepackage[nameinlink,noabbrev]{cleveref}
\usepackage{autonum}

\newtheorem{theorem}{Theorem}[section]

\newtheorem{proposition}[theorem]{Proposition}
\newtheorem{definition}[theorem]{Definition}
\newtheorem{corollary}[theorem]{Corollary}
\newtheorem{lemma}[theorem]{Lemma}

\newtheorem{ithm}{Theorem}[section]

\DeclareMathOperator{\End}{End}

\def\cal#1{\mathcal{#1}}
\def\bb#1{\mathbb{#1}}

\def\lr#1{\left\langle #1\right\rangle}

\newcommand{\h}{\frac{1}{2}}

\newcommand{\ga}{\textsl{g}}

\title{Fatness and positive curvatures on Riemannian submersions}

\author[Cavenaghi]{Leonardo F. Cavenaghi}
\address{Institute of Mathematics and Informatics, Bulgarian Academy of Sciences
 \& College of Arts and Sciences, Department of Mathematics, University of Miami, Ungar Bldg, 1365 Memorial Dr 515, Coral Gables, FL 33146, USA}
\email{leonardofcavenaghi@gmail.com}

\author[Grama]{Lino Grama}
\address{Instituto de Matemática, Estatística e Computação Científica (IMECC) da Universidade Estadual de Campinas (Unicamp), Cidade Universitária, Campinas - SP, 13083-856, Brazil} 
\email{linograma@gmail.com}

\author[Sperança]{Llohann D. Sperança}
\address{ICT - Unifesp, São José dos Campos, SP, Brazil}
\email{lsperanca@gmail.com}

\begin{document}

\maketitle



\begin{center}
\begin{minipage}{0.9\textwidth}
\small
\textbf{Abstract.}
Guided by rather counterintuitive properties of \emph{dual holonomy fields}, we introduce different notions of \emph{fatness} for Riemannian submersions. We introduce \emph{integral fatness} as an intrinsic condition implied by positive sectional curvature; \emph{sub-fatness} as a rigidity constraint that naturally arises in \emph{canonical variation} limits; and then develop the concept of \emph{weakly nonnegative curvature} (WNN), an intrinsic condition depending only on the horizontal distribution and the base metric. Positive sectional curvature implies fatness in this class, verifying Petersen--Wilhelm's Conjecture for homogeneous submersions and beyond. Finally, we introduce the notion of Riemannian submersion which are \emph{$A$-Clifford}: we show that them, like fatness and WNN, is intrinsic to the pair $(\mathcal H,b)$; that it singles out a canonical vertical gauge; and that in rank three it characterizes $3$-Sasakian geometry, yielding rigidity on $\mathrm{S}^3,\mathrm{SO}(3)$-principal bundles. 
\end{minipage}
\end{center}


\section{Introduction}

Riemannian submersions occupy a central place in Riemannian geometry, particularly as a tool for constructing and analyzing manifolds with lower curvature bounds. By the Horizontal Curvature Equation \cite{oneill,gray1967pseudo}, any lower sectional curvature bound on the total space of a Riemannian submersion \(\pi\colon M\to B\) descends to its base \(B\). 

With their natural metrics, all known closed manifolds with positive sectional curvature arise as the \emph{base space} of a Riemannian submersion \cite{berger1961varietes,wallach1972compact,eschenburg1982new,bazaikin1995one,Eschenburg-Heintze,Zillerpositive,escher2014topology}, except for the $7$-dimensional cohomogeneity-one example from \cite{dearricott20117} (and possibly one further sporadic example \cite{pw}). Nevertheless, positively curved manifolds remain exceedingly rare, and a major open challenge is to distinguish those simply connected manifolds admitting nonnegative sectional curvature from those supporting strictly positive curvature. Aiming for a better understanding of the possible ``scarcity'' or rigidity of positively curved Riemannian metrics, a striking conjecture of Petersen and Wilhelm is as follows:
If \(F\hookrightarrow M\xrightarrow{\,\pi\,}B\) is a Riemannian submersion from a closed manifold \(M\) with \(\sec(M)>0\), then $\dim F<\dim B.$

Although \emph{fatness} \cite{weinstein1980fat} is known to imply the Petersen–Wilhelm conjecture, and a lot more, little is known between its relation with the aforementioned underlying motivation: existence of metrics with positive sectional curvature. For example, the emergence of fatness under curvature assumptions has remained elusive. In this paper, we work on both ways: we show that fatness appears as a necessary condition under some assumption; and provide a weaker notion of fatness that is implied by positive sectional curvature.

\vspace{1em}

In a standard setting, one starts with a Riemannian submersion $\pi:(M,\ga)\to(B,b)$ between two Riemannian manifolds $(M,\ga),~(B,b)$ and then defines its \emph{horizontal distribution} (or \emph{connection}) as $\mathcal H = \mathcal{V}^{\perp}$, where $\mathcal V =\ker \mathrm d\pi$ is the \emph{vertical space}. In this work, we emphasize on conditions that do not depend on the entire Riemannian metric $\ga$. Therefore, we invite the reader to consider a horizontal distribution $\mathcal H$ over a submersion $M\to (B,b)$ as a smooth complement of $\mathcal V=\ker\mathrm d\pi$ and call $\ga$ \emph{adapted} to $(\mathcal{H}, b)$ if $\pi:(M,\ga)\to(B,b)$ defines a Riemannian submersion. That is, $\ga$ is \emph{adapted} to $(\mathcal{H}, b)$ if $\ga(\mathcal{V}, \mathcal{H}) = 0$ and $\ga(X,Y) = b(\mathrm d\pi X,\mathrm d\pi Y)$ for all $X\in\mathcal H$ and $Y\in\mathcal V$.
Throughout this paper, all manifolds are assumed to be compact and connected. 

\vspace{1em}

Important concepts in the literature do not rely on the submersion metric, but only on the horizontal connection. Arguably, one of the most important of those concepts is the above-mentioned \emph{fatness}: a horizontal distribution $\mathcal H$ is \emph{fat} if, at every $p\in M$ and every nonzero $X\in\mathcal H_p$, one has $\mathcal V_p=[\widetilde X,\mathcal H]^{v}$ for some horizontal extension $\widetilde X$ of $X$. 

A convenient alternative way to define fatness is through Gray--O'Neill's tensor. Given a horizontal vector $X\in\mathcal H$, the Gray--O'Neill tensor $A_X$ is a linear map $\mathcal H\to\mathcal V$ defined by $A_X Y = \h\big[ \widetilde X, \widetilde Y \big]^{v}$, where $\widetilde X$ and $\widetilde Y$ are horizontal extensions of $X$ and $Y$, respectively. Given a metric $\ga$ adapted to $(\mathcal{H}, b)$, one can consider the dual $A^*_X:\mathcal V\to\mathcal H$. Then, a horizontal distribution $\mathcal H$ is fat if, and only if, $A^*_X V \neq 0$ for all $X\in\mathcal H$ and $V\in\mathcal V$. Noticing that $X\mapsto A^*_XV$ is skew-symmetric, we conclude that fatness implies the Petersen--Wilhelm conjecture (see, e.g.,~\cite{amann2015positive,chen2016riemannian,gonzalez2016soft,jimenez2005riemannian}), but whether fatness emerges under curvature assumptions has remained an open question.

\vspace{1em}

This paper develops \emph{weaker} notions of fatness that are detected by vertizontal curvature. Section~\ref{sec:integral-fatness} introduces \emph{integral fatness}, a time-averaged coercivity condition on $|A^*_X u|^2$ along dual holonomy fields. We prove that it is implied by positive vertizontal curvature and that it implies a \emph{single dual leaf} \cite{wilkilng-dual}. 

Section~\ref{sec:subfat} derives a curvature condition naturally arising from \emph{canonical variation}, which we call \emph{sub-fatness}. This condition can be interpreted as a surviving non-negative sectional curvature condition on ghost vertizontal planes after fibers collapse to points in a canonical variation. We show, via a Riccati equation, that sub-fatness implies totally geodesic fibers. 

Section~\ref{sec:wnn} then develops weakly nonnegative curvature (WNN), a curvature condition that does not depend on the vertical metric. I.e., it is an intrinsic condition on the pair $(\mathcal H,b)$. We prove that positive sectional curvature implies fatness under this assumption. Notably, this assumption is satisfied whenever $(\mathcal H,b)$ admits a metric with totally geodesic fibers with non-negative sectional curvature, as in the case of cosets in normal homogeneous spaces.

Finally, Section~\ref{sec:aclifford} studies an extreme case of weakly non-negative submersions: \emph{$A$-Clifford} Riemannian submersions.

\vspace{1em}

Our main results are as follows and depend on the definition of a dual holonomy field. Let $c$ be a horizontal curve, i.e.\ $\dot c\in\mathcal H_{c(t)}$. A \emph{holonomy field} $\xi$ along $c$ is a vertical field satisfying
\begin{equation}\label{eq:holonomyfield}
\nabla_{\dot c}\xi=-A^*_{\dot c}\xi-S_{\dot c}\xi.
\end{equation}
A \emph{dual holonomy field} $\nu$ along $c$ satisfies
\begin{equation}\label{eq:dualhol}
\nabla_{\dot c}\nu=-A^*_{\dot c}\nu+S_{\dot c}\nu.
\end{equation}
Here, $S_{\dot c}:\cal V\rightarrow \cal V$, $S_{\dot c}\xi = -(\nabla_{\xi} \widetilde{X})^v$, denotes the second fundamental form operator where $\widetilde{X}$ is an arbitrary horizontal extension of $\dot c$. Both equations are linear with smooth coefficients along a complete horizontal curve, so solutions exist for all time and are unique; in particular, a solution with non-zero initial condition never vanishes.

\begin{ithm}\label{ithm:dual-holonomy-fields}
  Let $\pi:(M,\mathcal H)\to (B,b)$ be a Riemannian submersion with adapted metric $\ga$. Let $\xi,\nu$ be a holonomy field and a dual holonomy field, respectively, along a horizontal curve $c$, $\dot c=X$. Then,
  \begin{enumerate}
    \item $X\ga(\xi,\nu)=0$
    \item Let $\ga'$ be  adapted to $(\mathcal H,b)$ and denote its dual Gray--O'Neill-tensor as $A^\dagger$. Then there is $\nu'$, a dual holonomy field with respect to $\ga'$ along $c$, such that 
    \[ A^*_X\nu = A^\dagger_X\nu'\]
    along $c$. 
    \item\label{item:log-curvature} Let $u=\nu/|\nu|$, and $z=\log|\nu|$. Then $z'=\langle S_{\dot\gamma}u,u\rangle$ and
    \begin{equation}\label{eq:log-curvature-intro}
    z''=K(\dot\gamma,u)+3|S_{\dot\gamma}u|^2-2\langle S_{\dot\gamma}u,u\rangle^2-|A^*_{\dot\gamma}u|^2.
    \end{equation}
    In particular,
    \[z'' \geq  (z')^2+  K(\dot\gamma,u)-|A^*_{\dot\gamma}u|^2.\]
    
  \end{enumerate}  
\end{ithm}

\begin{ithm}\label{ithm:integral-fatness-intro}
For a submersion $\pi: (M, \mathcal{H}) \rightarrow (B, b)$ with an adapted metric $\ga$, we say that $\ga$ is \emph{integrally fat} if 
\[\liminf_{T\to\infty}\frac1T\int_0^T |A^*_{\dot c} u|^2\,\mathrm dt > 0\]
for all $X\in\mathcal H$ and $u=\nu/|\nu|$, for $\nu\in\mathcal V$ a dual holonomy field along the geodesic defined by $\dot c = X$.
\begin{enumerate}
\item $\ga$ is integrally fat if and only if every metric adapted to $(\cal H,b)$ is integrally fat. \label{item:independent-vert}
\item If $\ga$ has positive vertizontal curvature, then $\pi$ is integrally fat.
\item If $\pi$ is integrally fat, then $\pi$ has a single dual leaf. I.e., every two points can be joined by a continuous, piecewise $C^1$ horizontal curve.
\end{enumerate}
\end{ithm}

In particular, if $\ga$ has positive vertizontal curvature, then $\pi$ has a single dual leaf. This strengthens \cite{speranca2011anote} by removing the bounded holonomy hypothesis.

Second, we define the sub-fatness regime which arises from canonical variation. We say that $\pi$ is \emph{sub-fat} if 
\[ K(X,u) \geq |A^*_X u|^2\]
for all $X\in\mathcal H$ and $u\in\mathcal V$ where $K$ is the unreduced sectional curvature $K(X,u)=R_\ga(X,u,u,X)$.
\begin{ithm}\label{ithm:subfat-rigidity-intro}
If $\pi$ is sub-fat, then its fibers are totally geodesic and $K(X,u)=|A^*_X u|^2$ for all vertizontal pairs.
\end{ithm}

The beauty of Theorem \ref{ithm:subfat-rigidity-intro} is in its proof which relies on a Riccati equation that arises from the dual holonomy field curvature equation.
Third, we introduce weakly nonnegative curvature (WNN), a curvature condition that does not depend on the vertical metric. I.e., it is an intrinsic condition on the pair $(\mathcal H,b)$. We prove that positive sectional curvature implies fatness under this assumption. Notably, this assumption is satisfied whenever $(\mathcal H,b)$ admits a metric with totally geodesic fibers with non-negative sectional curvature, as in the case of cosets in normal homogeneous spaces.

\begin{definition}\label{dfn:wnn}
An adapted metric $\ga$ to $(\mathcal H,b)$ is called \emph{weakly nonnegatively curved} (WNN) if every $p\in M$ has a neighborhood $U$ and a constant $\tau>0$ such that
\begin{equation}\label{eq:WNN-intro}
\tau\, |X|\,|A^*_X\xi|^2 \geq \langle (\nabla_X A^*)_X\xi + A^*_X S_X\xi, A^*_X\xi \rangle
\end{equation}
for all $X \in \mathcal H|_U$, $\xi \in \mathcal V|_U$. A pair $(\mathcal H,b)$ is \emph{weakly nonnegatively curved} if it admits a WNN adapted metric.
\end{definition}

The interesting point in this inequality arises when one takes $\xi$ as a dual holonomy field along a geodesic $c$, $\dot c =X$. In this case,
\begin{equation}\label{eq:mixed-curvature-is-a-norm}\langle (\nabla_X A^*)_X\xi + A^*_X S_X\xi, A^*_X\xi \rangle = \frac{1}{2}|A^*_X\xi|^2,\end{equation}
making Gronwall's inequality a natural tool. Moreover, this term is inherent to $(\cal H, b)$, i.e., it is independent of the vertical metric.

Both sides of \eqref{eq:WNN-intro} are homogeneous of degree three in $X$, so it suffices to require the inequality for unit horizontal vectors, in which case it reads $\tau|A^*_X\xi|^2\geq \lr{(\nabla_X A^*)_X\xi + A^*_X S_X\xi, A^*_X\xi}$.

\begin{ithm}\label{thm:varintro}
Let $\pi: (M, \mathcal{H}) \rightarrow (B, b)$ be  a submersion with an adapted metric $\ga$. Then:
\begin{enumerate}
 	\item 	$\ga$ is weakly nonnegatively curved if and only if every metric adapted to $(\cal H,b)$ is weakly nonnegatively curved \label{item:1}
	\item  if $\ga$ is nonnegatively curved and the fibers of $\pi$ are totally geodesic fibers, then $\ga$ is weakly nonnegatively curved \label{item:2}
	\item Suppose $(M,\ga)$ is complete, weakly nonnegatively curved and has positive sectional curvature. Then $\cal H$ is fat. \label{item:essevainofim}
\end{enumerate}
\end{ithm}

As an immediate consequence, we recover the following result from the literature:
\begin{ithm}\label{ithm:homo}
Let $K < H < G$ be compact Lie groups and $\pi: (G/K, \mathcal{H}) \rightarrow (G/H, b)$ the right coset projection. Let $b$ be the normal homogeneous metric in $G/H$ and $\mathcal H$ as the standard $G$-invariant horizontal distribution. If there is a metric of positive sectional curvature  adapted to $(\mathcal H, b)$, then $\mathcal{H}$ is fat.
\end{ithm}
Gonz\'alez and Radeschi~\cite{gonzalez2017note} studied submersions from spaces homotopy equivalent to the known positively curved examples without curvature assumptions. In contrast, Theorem \ref{ithm:homo} imposes a curvature condition but may apply to different spaces; for example, with $G = \mathrm{Spin}(5)$, $H = \mathrm{Spin}(3)$, and the Berger space $B^7 = G/H$~\cite{berger1961varietes}, Theorem \ref{ithm:homo} implies that $\mathrm{Spin}(5)$ cannot admit a positively curved metric adapted to $(\mathcal{H}, b)$. Theorem~\ref{thm:varintro} encompasses a broader class, including the non-standard examples of~\cite{kerin-shankar}. Moreover, the present proof uses no topological or Lie-theoretic classification.

We recall that fatness alone imposes restrictions on dimensions and characteristic numbers. By combining Theorem~\ref{thm:varintro}\ref{item:essevainofim} with the obstructions in Florit and Ziller \cite{Florit2011}, we conclude that WNN is one of the few extra curvature-related conditions that transform invariants into obstructions to the existence of positive sectional curvature. In particular, a WNN submersion with positive sectional curvature satisfies the Radon--Hurwitz bound $\dim F\le\rho(\dim B)-1$, 
and all its Weinstein invariants are non-zero \cite{Florit2011}. For instance, by \cite[Corollaries 2--4]{Florit2011}, the unit tangent bundles of $\mathrm S^8$, $\bb{CP}^4$ and $\bb{CP}^2$ carry no submersion metric with positive sectional curvature adapted to  their natural $(\cal H, b)$.

\vspace{1em}

Finally, we introduce the concept of $A$-Clifford Riemannian submersions and apply them to principal bundles. Recall from \cite{MOROIANU2011940} that a \emph{rank $r$ Clifford structure} on a Riemannian manifold $(B,b)$ is an oriented rank $r$ Euclidean bundle $(E,\ga_E)$ over $B$, together with a non-vanishing algebra bundle morphism $\varphi:\mathrm{Cl}(E,\ga_E)\to\mathrm{End}(TB)$ into skew-symmetric endomorphisms. When $\pi:(M,\ga)\to(B,b)$ is a Riemannian submersion, the Gray--O'Neill tensor defines skew-symmetric endomorphisms: given $u\in\cal V$, let $A^u:\cal H\to\cal H$ be $A^uX=A^*_Xu$. We study the class of Riemannian submersions for which $u\mapsto A^u$ defines a Clifford morphism.

Namely, $\pi$ is \emph{$A$-Clifford} if $(A^u)^2$ is a multiple of the identity for every $u\in\cal V$, say $(A^u)^2=-N(u)\mathrm{Id}$. We call it \emph{$A$-Clifford in the normalized gauge} if $N(u)=|u|^2$. The latter accounts exactly for the existence of an algebra bundle morphism of $\mathrm{Cl}(\cal V,\ga|_{\cal V})$:  
\[u\cdot u=-\ga(u,u)1\] 
The distinction matters: the first condition is intrinsic to $(\cal H,b)$, while the second fixes a canonical vertical metric.

\begin{ithm}\label{ithm:aclifford}
  Let $\pi:(M,\ga)\rightarrow (B,b)$ be a Riemannian submersion. Then $\ga$ is $A$-Clifford if and only if it is $A$-Clifford for every metric adapted to $(\cal H,b)$; the condition depends only on $(\cal H,b)$. Moreover:
  \begin{enumerate}
    \item $N$ is a quadratic form on $\cal V$ and $|A^*_Xu|^2=N(u)|X|^2$; hence $\pi$ is fat if and only if $N$ is definite;
    \item if $\pi$ is fat, then $N$ is the unique adapted vertical metric in which the structure is normalized, and in that gauge $u\mapsto A^u$ is a rank $r$ Clifford structure with $K(X,u)=|u|^2|X|^2$;
    \item if $\pi$ is fat, then $\ga$ is weakly nonnegatively curved; 
    \item in the normalized gauge one may take $\tau$ as the maximum operator norm of $S_X$ for unit $X\in \mathcal{H}$. In particular, the right-hand side of \eqref{eq:WNN-intro} vanishes identically for totally geodesic fibers, and $N(\nu)=|\nu|^2$ is constant on holonomy fields.
  \end{enumerate}
\end{ithm}

The $A$-Clifford metrics thus lie on the extreme side of weak nonnegative curvature: in the normalized gauge with totally geodesic fibers, the right-hand side of \eqref{eq:WNN-intro} is zero. Behind items (3)-(4), and the rigidity results of Section~\ref{sec:aclifford}, is the following pair of curvature identities under totally geodesic fibers (\cite{gw} page 44):
 \begin{equation}
\begin{aligned}
    R(X,\xi,A^*_X\xi,Y) &= \langle(\nabla_XA^*)_Y\xi,A^*_X\xi\rangle.
\end{aligned}
\end{equation}
In particular:
 \begin{equation}
\begin{aligned}
R(X,\xi,A^*_X\xi,X)=\dfrac{1}{2}\dfrac{\mathrm d}{\mathrm dt}\Big|_{t=0}|A^*_{\dot c}\xi(t)|^2,
\end{aligned}
\end{equation}
along the unit speed geodesic $c$ generated by $X$, and a holonomy field $\xi(t)$ extending $\xi=\xi(0)$.

One notices that the normalized $A$-Clifford condition generalizes $3$-Sasakian structures. In fact, they meet in rank three.

\begin{ithm}\label{ithm:3sas}
  Let $\pi:G\hookrightarrow (P,\ga)\to (B,b)$ be a Riemannian submersion which is $A$-Clifford in the normalized gauge where $G$ acts by isometries. Suppose $G=\mathrm S^3$ or $\mathrm{SO}(3)$ and that $\pi$ has totally geodesic fibers of curvature one. Then $(P,\ga)$ is $3$-Sasakian if and only if the associated $\mathrm{Cl}_{0,3}$-module is isotypic. In particular $(B,b)$ is quaternionic K\"ahler with positive scalar curvature. 
\end{ithm}

The converse part of Theorem \ref{ithm:3sas} directly relates to \cite{Kashiwada2001}: an \textit{isotypic} $\mathrm{Cl}_{0,3}$-module \cite{lawson2016spin} guarantees the existence of a \textit{contact 3-structure}, which, by \cite{Kashiwada2001}, is 3-Sasakian.

\vspace{1em}

\section{Dual holonomy fields and integral fatness}
\label{sec:integral-fatness}

Let $\pi:(M,\ga)\to(B,b)$ be a Riemannian submersion with vertical bundle $\mathcal{V}\subset TM$ and horizontal distribution $\mathcal H$. We write $K(X,u)=R_\ga(X,u,u,X)$ for the \emph{unreduced} sectional curvature of a vertizontal plane spanned by horizontal $X$ and vertical $u$. The Gray--O'Neill integrability tensor is
\[
A_X Y = \h\big[ \widetilde X, \widetilde Y \big]^{v} = (\nabla_{ X} \widetilde Y )^v = (\nabla_{ Y} \widetilde X)^v,
\]
for horizontal extensions $\widetilde X,\widetilde Y$ of $X,Y\in\mathcal H$, and $A^*$ is its metric adjoint. We define the \emph{second fundamental form} of the fibers as
\[
S_X\xi = -(\nabla_\xi \widetilde X)^{v},
\]
for $X\in\mathcal H$, $\xi\in\mathcal V$. We say that $\pi$ has \emph{totally geodesic fibers} if $S\equiv 0$.

By compactness of $M$, set
\[
\begin{aligned}
    \sigma   &:= \max_{X\in\mathcal H,\,|X|=1}|S_X|_{\mathrm{op}} < \infty, \\[1ex]
    \kappa_v &:= \min\big\{K(X,u): X\in\mathcal H_p,\,u\in\mathcal V_p,\,|X|=|u|=1,\,p\in M\big\}.
\end{aligned}
\]
Here $|\cdot|_\mathrm{op}$ is the operator norm. We see that $\pi$ has \emph{positive vertizontal curvature} when $\kappa_v>0$.

\subsection{Holonomy and dual holonomy fields}

Given a horizontal geodesic $c$, its \emph{infinitesimal holonomy transformation} $\hat c(t):\mathcal V_{c(0)}\to\mathcal V_{c(t)}$ is defined by $\hat c(t)\xi(0)=\xi(t)$ for the holonomy field $\xi$ with $\xi(0)$ prescribed. Dual holonomy fields satisfy $\nu(t)=\hat c(t)^{-\ast}\nu(0)$, where $\hat c(t)^{-*}$ is the inverse of the metric dual of $\hat c(t)$, $\hat c(t)^{-*} = (\hat c(t)^{-1})^*$ \cite[Proposition~4.1]{speranca2017on}.

For a horizontal geodesic $\gamma:\mathbb R\to M$, define
\[
C(\gamma):=\mathrm{span}\big\{\hat\gamma(t)^{-1}(A_{\dot\gamma(t)}Z): Z\in\mathcal H_{\gamma(t)},\,t\in\mathbb R\big\}\subset\mathcal V_{\gamma(0)}.
\]
The \emph{dual leaves} \cite{wilkilng-dual,speranca2017on} are the immersed submanifolds tangent to the horizontal distribution and to $\mathcal V\cap T L^\#$; a kind of first-order reachability along horizontal curves is measured by the sets $C(\gamma)$ (along the lines of \cite{ambrose-singer}).

Let $\ga,\ga'$ be metrics adapted to $(\mathcal H,b)$ and let $P_p:\mathcal V_p\to\mathcal V_p$ be defined by $\ga'(\xi,\eta)=\ga(P_p\xi,\eta)$. We use $^*$ for $\ga$-duals and $^\dagger$ for $\ga'$-duals.

\begin{lemma}\label{lem:dualinv}
	If $\nu'(t)$ is a dual holonomy field on $(M,\ga')$ along a horizontal curve $c$, then $\nu'(t)=P^{-1}_{c(t)}\nu(t)$, where $\nu(t)$ is the dual holonomy field on $(M,\ga)$ with $\nu(0)=P_{c(0)}\nu'(0)$. In particular,
	\begin{equation}\label{eq:lemdualinv}
	A^\dagger_{\dot c}\nu'(t)=A^*_{\dot c}\nu(t).
	\end{equation}
\end{lemma}
\begin{proof} Note that holonomy transport is completely determined by $\cal H$. Therefore, so are holonomy fields. Moreover, it is sufficient to investigate how inifinitesimal holonomy transformations trasform under the change of metric. Given $u,v\in\cal V$
\begin{align}
    \ga'(\hat c(t)^{-\dagger}u, v) &= \ga'(u, \hat c(t)^{-1}v) = \ga(Pu, \hat c(t)^{-1}v) \\
    &= \ga(\hat c(t)^{-\ast}Pu, v) = \ga'(P^{-1}\hat c(t)^{-\ast}Pu, v).
\end{align}
In particular,
\[
P^{-1}_{c(t)}\hat c(t)^{-\ast}P_{c(0)}\nu'(0)=P^{-1}_{c(t)}\nu(t).\] 
On the other hand,  $A^\dagger_X\nu'=A^*_X P_p\nu'=A^*_X \nu$.
\end{proof}

\subsection{The logarithmic curvature identity}

We now derive the main equations that relate dual holonomy fields and curvature (Theorem \ref{ithm:dual-holonomy-fields}\ref{item:log-curvature}).

\begin{lemma}\label{lem:norm-derivative}
Let $\nu$ be a dual holonomy field along a horizontal curve $c$. Then,
\[
\frac{\mathrm{d}}{\mathrm{d}t}|\nu|^2 = 2\langle S_{\dot c}\nu,\nu\rangle.
\]
In particular, along a unit speed horizontal curve, $\big|(\log|\nu|)'\big|\le\sigma$ wherever $\nu\neq 0$.
\end{lemma}
\begin{proof}
$\frac{\mathrm{d}}{\mathrm{d}t}|\nu|^2=2\langle\nabla_{\dot c}\nu,\nu\rangle=2\langle -A^*_{\dot c}\nu+S_{\dot c}\nu,\nu\rangle$, and $A^*_{\dot c}\nu\in\mathcal H$ is orthogonal to $\nu$.
\end{proof}

\begin{lemma}\label{lem:log-curvature}
Let $\gamma$ be a unit-speed horizontal geodesic, $\nu\not\equiv 0$ a dual holonomy field along $\gamma$, $u=\nu/|\nu|$, and $g=\log|\nu|$. Then
\begin{equation}\label{eq:log-curvature}
g'=\langle S_{\dot\gamma}u,u\rangle,\qquad
g''=K(\dot\gamma,u)+3|S_{\dot\gamma}u|^2-2\langle S_{\dot\gamma}u,u\rangle^2-|A^*_{\dot\gamma}u|^2.
\end{equation}
Consequently,
\begin{equation}\label{eq:log-curvature-integrated}
g''\ge K(\dot\gamma,u)+(g')^2-|A^*_{\dot\gamma}u|^2\ge \kappa_v+(g')^2-|A^*_{\dot\gamma}u|^2.
\end{equation}
\end{lemma}
\begin{proof}
Write $f=|\nu|^2=e^{2g}$. Then $f'=2g'e^{2g}=2g'f$ and $f''=(2g''+4(g')^2)f$, so
\[
\tfrac12 f''=\bigl(g''+2(g')^2\bigr)f .
\]
By \ref{lem:norm-derivative}, $f'=2\langle\nabla_{\dot\gamma}\nu,\nu\rangle=2\langle S_{\dot\gamma}\nu,\nu\rangle$, so
\[
g'=\frac{f'}{2f}=\frac{\langle S_{\dot\gamma}\nu,\nu\rangle}{|\nu|^2}=\langle S_{\dot\gamma}u,u\rangle .
\]
\cite[Proposition~4.2]{speranca2017on} gives, for any dual holonomy field,
\begin{equation}\label{eq:K}
K_\ga(\dot\gamma,\nu)=\tfrac12 f''-3|S_{\dot\gamma}\nu|^2+|A^*_{\dot\gamma}\nu|^2.
\end{equation}
Dividing \eqref{eq:K} by $f=|\nu|^2$ and using that $K_\ga(\dot\gamma,\nu)$, $|S_{\dot\gamma}\nu|^2$, $|A^*_{\dot\gamma}\nu|^2$ are all quadratic in $\nu$,
\[
K(\dot\gamma,u)=\dfrac{1}{2}\frac{ f''}{f}-3|S_{\dot\gamma}u|^2+|A^*_{\dot\gamma}u|^2
=g''+2(g')^2-3|S_{\dot\gamma}u|^2+|A^*_{\dot\gamma}u|^2 ,
\]
where the second equality substitutes $\tfrac12 f''=(g''+2(g')^2)f$. Then,
\[
g''=K(\dot\gamma,u)-2\langle S_{\dot\gamma}u,u\rangle^2+3|S_{\dot\gamma}u|^2-|A^*_{\dot\gamma}u|^2 ,
\]
which is \eqref{eq:log-curvature}. For \eqref{eq:log-curvature-integrated}, Cauchy--Schwarz with $|u|=1$ gives $\langle S_{\dot\gamma}u,u\rangle^2\le|S_{\dot\gamma}u|^2$, so
\[
3|S_{\dot\gamma}u|^2-2\langle S_{\dot\gamma}u,u\rangle^2\ge \langle S_{\dot\gamma}u,u\rangle^2=(g')^2,
\]
whence $g''\ge K(\dot\gamma,u)+(g')^2-|A^*_{\dot\gamma}u|^2$. Finally $K(\dot\gamma,u)\ge\kappa_v$ since $\dot\gamma,u$ are orthonormal.
\end{proof}


\subsection{Integral fatness}

\begin{definition}\label{dfn:integral-fatness}
A Riemannian submersion is \emph{integrally fat} if, for every unit-speed horizontal geodesic $\gamma$ and every non-zero dual holonomy field $\nu$ along $\gamma$ with unit direction $u(t)=\nu(t)/|\nu(t)|$,
\[
\liminf_{T\to\infty}\frac{1}{T}\int_0^T |A^*_{\dot\gamma(t)}u(t)|^2\,\mathrm{d}t > 0.
\]
\end{definition}

Observe that, by Lemma \ref{lem:dualinv}, the integrand is independent of the choice of vertical metric. More specifically, following the conventions in its statement,
\[\int_0^T |A^*_{\dot\gamma(t)}u(t)|_{\ga}^2\,\mathrm{d}t = \int_0^T |A^{\dagger}_{\dot\gamma(t)}u'(t)|_{\ga'}^2\,\mathrm{d}t,\]
thus proving Theorem \ref{ithm:integral-fatness-intro}\ref{item:independent-vert}.

\begin{theorem}\label{thm:integral-fatness-estimate}
If $\kappa_v>0$, then for every unit speed horizontal geodesic $\gamma$, every non-zero dual holonomy field $\nu$ along $\gamma$, and every $T>0$,
\[
\int_0^T |A^*_{\dot\gamma(t)}u(t)|^2\,\mathrm{d}t \ge \kappa_v T - 2\sigma,
\qquad u=\nu/|\nu|.
\]
In particular, $\pi$ is integrally fat.
\end{theorem}
\begin{proof}
By \eqref{eq:log-curvature-integrated}, at every $t\in[0,T]$,
\[
|A^*_{\dot\gamma}u|^2\ge \kappa_v+(g')^2-g''\ge \kappa_v-g'' .
\]
Integrating over $[0,T]$,
\[
\int_0^T |A^*_{\dot\gamma}u|^2\,\mathrm dt\ \ge\ \kappa_v T-\int_0^T g''\,\mathrm dt
\ =\ \kappa_v T-\bigl(g'(T)-g'(0)\bigr).
\]
By Lemma \ref{lem:norm-derivative}, $|g'|\le\sigma$ pointwise, so $g'(T)-g'(0)\le 2\sigma$, giving
\[
\int_0^T |A^*_{\dot\gamma}u|^2\,\mathrm dt\ \ge\ \kappa_v T-2\sigma .
\]
Whence, $\liminf_{T\to\infty}\frac1T\int_0^T|A^*_{\dot\gamma}u|^2\,\mathrm dt\ge\kappa_v>0$.
\end{proof}

\begin{corollary}\label{cor:C-equals-V}
If $\pi$ is integrally fat, then $C(\gamma)=\mathcal V_{\gamma(0)}$ for every horizontal geodesic $\gamma$.
\end{corollary}
\begin{proof}
Let $\nu$ be a dual holonomy field along $\gamma$ with $\nu(0)= \nu_0\perp C(\gamma)$. From \cite[Lemma~6.4]{speranca2017on}, we know that $A^*_{\dot\gamma(t)}\nu(t)=0$ for all $t$.
Hence, if $0\neq\nu_0\perp C(\gamma)$, $\liminf_{T\to\infty}\frac{1}{T}\int_0^T|A^*_{\dot\gamma}u|^2\,\mathrm{d}t=0$, contradicting Definition \ref{dfn:integral-fatness}.
\end{proof}

The proof of Theorem \ref{ithm:integral-fatness-intro} is concluded with the following result.

\begin{theorem}\label{thm:single-dual-leaf}
If $\pi$ is integrally fat, then $\pi$ has a \emph{single dual leaf}.
\end{theorem}
\begin{proof}
Fix $q\in M$. Let $\mathcal A_q=\mathrm{span}\{h^{-1}(A_XY): X,Y\in\mathcal H_{\tau(h)},\,h\in\mathcal E_q\}$, where $\mathcal E_q$ is the holonomy group at $q$ \cite{speranca2017on}. By Corollary \ref{cor:C-equals-V}, $\mathcal V_q=C(\gamma)\subseteq\mathcal A_q$ for every horizontal geodesic through $q$. Horizontal fields are tangent to $L^\#$ \cite{wilkilng-dual}, and $2A_XY=[X,Y]^{v}$ is tangent to dual leaves; infinitesimal holonomy transport preserves $\mathcal V\cap TL^\#$ \cite[Section~3]{speranca2017on}. Thus $T_qL^\#\supseteq\mathcal H_q\oplus\mathcal V_q=T_qM$, so each dual leaf is open in the compact manifold $M$ and there is at least one.
\end{proof}

\section{Sub-fatness and canonical variation}
\label{sec:subfat}

Canonical variation rescales the vertical metric while fixing the horizontal distribution. Its vertizontal curvature limit motivates a pointwise  condition that we call \emph{sub-fatness}. It rigidifies the submersion through the same Riccati identity of Section~\ref{sec:integral-fatness}.

\subsection{Canonical variation}

Let $\ga$ be  adapted to $(\cal H, b)$ and $\mu>0$. The \emph{canonical variation} $\ga_\mu$ is defined by
\[
\ga_\mu|_{\mathcal H}=\ga|_{\mathcal H},\qquad \ga_\mu|_{\mathcal V}=\mu\,\ga|_{\mathcal V},\qquad \ga_\mu(\mathcal H,\mathcal V)=0.
\]
The Koszul formula on horizontal and vertical vectors shows that $A^{\ga_\mu}=A$, $S^{\ga_\mu}=S$, and $(A^*)^{\ga_\mu}=\mu A^*$  (\cite[Section 2.1]{gw}, where one takes $e^{2t}=\mu$ there, with a constant $t$.)

\begin{proposition}\label{prop:canonical-variation}
For $\ga_\mu$ as above, $\ga_\mu$ is adapted to $(\cal H,b)$ and:
\begin{enumerate}
\item \label{item:dual-holonomy-coincide} The sets of  dual holonomy fields in $\ga$ and $\ga_\mu$ coincide. Namely, $\nu$ is a dual holonomy field in $\ga$ if and only it is a dual holonomy field in $\ga'$.
\item\label{item:canonical-K} The unreduced vertizontal sectional curvature transforms as
\[
K_{\ga_\mu}(X,T)=\mu K_\ga(X,T)+\mu(\mu-1)|A^*_X T|_{\ga}^2
\]
for a $\ga$-orthonormal vertizontal pair $(X,T)$.
\item \label{item:canonical-sec-limit} As $\mu\to 0^+$, the normalized vertizontal curvature $\mathrm{sec}_{\ga_\mu}(X\wedge T)$ satisfies 
\[\lim_{\mu\to 0^+}\mathrm{sec}_{\ga_\mu}(X\wedge T) = K_\ga(X,T)-|A^*_X T|_{\ga}^2.\]
\end{enumerate}
\end{proposition}
\begin{proof}
Item \ref{item:dual-holonomy-coincide} follows from Lemma \ref{lem:dualinv}. Item \ref{item:canonical-K} is the vertizontal specialization of \cite[Example 2.1.1]{gw}, and \ref{item:canonical-sec-limit} follows by observing that $|u|^2_{\ga _{\mu}} = \mu|u|^2_{\ga}$.
\end{proof}

Since $(A^*)^{\ga_\mu}=\mu A^*$ and $|u|_{\ga_\mu}^2=\mu|u|_{\ga}^2$ for $u\in\mathcal V$, a unit vertizontal pair $(X,u)$ for $\ga$ gives
\begin{equation}\label{eq:subfat-canonical}
K_{\ga_\mu}(X,u)-|(A^*)^{\ga_\mu}_X u|_{\ga_\mu}^2
=\mu\bigl(K_\ga(X,u)-|A^*_X u|_{\ga}^2\bigr).
\end{equation}
In particular, $\ga_\mu$ is sub-fat if and only if $\ga$ is.

The collapse limit in Proposition~\ref{prop:canonical-variation}(3) motivates a positive/non-negative curvature condition on the surviving vertizontal curvature. That is, to require $K(X,u)\ge|A^*_X u|^2$. This is the exact geometric content of sub-fatness. In contrast with the other courvature conditions, we will see that sub-fatness is extremely sensitive to the vertical metric. In fact, it dictates totally geodesic fibers.

\subsection{Sub-fatness and rigidity}

\begin{definition}\label{dfn:subfat}
A Riemannian submersion is \emph{sub-fat} if
\begin{equation}\label{eq:subfat}
K(X,u)\ge |A^*_X u|^2
\end{equation}
for all $X\in\mathcal H$, $u\in\mathcal V$.
\end{definition}

Equality holds identically when the fibers are totally geodesic \cite[p.~28]{gw}. Sub-fatness is strictly stronger than the time-averaged condition of Definition \ref{dfn:integral-fatness} and, by \eqref{eq:subfat-canonical}, it is invariant under canonical variation.

\begin{theorem}\label{thm:subfat-rigidity}
If $\pi$ is sub-fat, then:
\begin{enumerate}
\item\label{item:subfat-tg} the fibers are totally geodesic ($S\equiv 0$);
\item\label{item:subfat-K-equals-A} $K(X,u)=|A^*_X u|^2$ for all vertizontal pairs.
\end{enumerate}
The same conclusions hold for every canonical variation $\ga_\mu$.
\end{theorem}
\begin{proof}
Fix $q\in M$ and a unit pair $(X,u)\in\mathcal H_q\times\mathcal V_q$ with respect to $\ga$. Let $\gamma(t)=\exp_q(tX)$ and let $\nu$ be the dual holonomy field with $\nu(0)=u$. Set $g=\log|\nu|_{\ga}$ and $z=g'$. 

Assume $\ga$ sub-fat. By the logarithmic curvature identity, and Lemma \ref{lem:norm-derivative}:
\[
z'\ge K_{\ga}(\dot\gamma,u(t))-|(A^*)^{\ga}_{\dot\gamma}u(t)|_{\ga}^2+z^2\ge z^2,
\qquad |y|\le\sigma
\]
If $z(t_0)>0$, then $z$ is non-decreasing and $(-1/z)'=z'/z^2\ge 1$ forces $z\to+\infty$ before time $t_0+1/z(t_0)$, contradicting $|z|\le\sigma$. Similarly $z(t_1)<0$ is impossible. Hence $z\equiv 0$, so $\langle S_X u,u\rangle=0$ for every unit pair and $S_X\equiv 0$. Then $z''\equiv 0$ and the logarithmic curvature identity for $\ga_\mu$ gives $K_{\ga}(X,u)=|(A^*)^{\ga}_X u|_{\ga}^2$.
\end{proof}

\vspace{1em}

	\section{An intrinsic curvature condition for submersions over Riemannian manifolds}
\label{sec:wnn}

In the Introduction, we have defined  weakly nonnegative curvature (WNN) (Definition~\ref{dfn:wnn}) for adapted metrics $\ga$. Like integral and sub-fatness, WNN is formulated through dual holonomy transport. Unlike the classical nonnegative sectional curvature in $M$, the WNN depends only on $(\mathcal H,b)$ and is strictly weaker than $\sec(M)\ge 0$ in the context of totally geodesic fibers, as it does not restrict the base curvature.

Below we prove Theorem \ref{thm:varintro}. Theorem \ref{ithm:homo} then follows as $G/K\to G/H$ has a metric as in Theorem \ref{thm:varintro}\textit{\ref{item:2}}. We start with the following auxiliary result.

\begin{lemma}\label{lem:flatgeo}
Let $\ga$ be a WNN metric adapted to $(\cal H,b)$. If $\nu$ is a dual holonomy field along a horizontal geodesic $c$ satisfying $A^*_{\dot c}\nu(0)=0$, then $A^*_{\dot c}\nu(t)=0$ for all $t$.
\end{lemma}
\begin{proof}
Take $|\dot c|=1$ and denote $u(t)=|A^*_{\dot c}\nu(t)|^2$. Equations \eqref{eq:mixed-curvature-is-a-norm} and \eqref{eq:WNN-intro} combined give $u'(t)\le 2\tau u(t)$; cf. \cite[Lemma~2.2]{speranca_WNN}. Gr\"onwall's inequality implies $u(t)\leq u(0)e^{2\tau t}$ for all $t\geq 0$. Therefore, if $u(0)=0$, $u(t)=0$ for all $t\geq0$. Replacing $c$ by $\tilde c(t)=c(-t)$ proves the assertion for all $t$.
\end{proof}

Moreover, if the vertizontal curvature of $M$ is non-negative, an argument as in the proof of Theorem \ref{thm:subfat-rigidity} shows that $S_X\nu\equiv 0$, thus $\nu$ is parallel.

\begin{proof}[Proof of item \ref{item:1}]
Let us show that the RHS of \ref{eq:WNN-intro} is invariant under changes of adapted metrics, up to a change in the dual holonomy field. More specifically, let $\ga,\ga'$ be adapted to $(\cal H,b)$ such that $\ga'(\xi,\eta)=\ga(P\xi,\eta)$. Denote the $\ga$- and $\ga'$-duals of $A$ as $A^*$ and $A^\dagger$. Moreover, let $\nabla^{\ga},~\nabla^{\ga}$ be the respective covariant derivatives and note that $(\nabla^\ga Z)^h = (\nabla^{\ga'} Z)^h$ for any horizontal field $Z\in\cal H$, since the metric in the base did not change. Then,
\[[\nabla^{\ga'}(A^\dagger_XP^{-1}\nu)]^h = [\nabla^{\ga}(A^\dagger_XP^{-1}\nu)]^h= [\nabla^{\ga}(A_X^*\nu)]^ h.\]

Nonetheless, for a dual holonomy field $\nu$ in $\ga$ along a horizontal geodesic $c=\exp(tX)$,
\[\nabla_X^\ga A^*_X\nu = (\nabla_X^\ga A^*)_X\nu + A^*_X(\nabla^\ga_X \nu)^v=(\nabla_X^\ga A^*)_X\nu+A^*_XS_X\nu.\]
Analogously, $A^\dagger_X P^{-1}\nu=A^*_X\nu$, concluding that the left-most term of the equation does not change. In particular, the WNN inequality holds with the same constant $\tau$. \qedhere
\end{proof}

Note that the proof above shows invariance under a relatively stronger inequality:
\[\tau|A^*_X\nu| \geq |(\nabla_X A^*)_X\nu+A^*_XS_X\nu|.\]
One might wonder if the new inequality is also inherited by non-negative sectional curvature. The answer is affirmative and both fact were used in \cite{speranca_grove} to conclude a result stronger than Lemma \ref{lem:flatgeo}.

\begin{proof}[Proof of item \ref{item:2}]
Item \ref{item:2} follows by analyzing the discriminant of the polynomial
\[K_{\ga}(X,\lambda Z+\xi)=K_{\ga}(\xi,X)+2\lambda R_{\ga}(X,Z,\xi,X)+\lambda^2K_{\ga}(X,Z),\]
with $Z=A^*_X\xi$. Since $(M,\ga)$ has nonnegative sectional curvature, this polynomial is everywhere positive and its discriminant is negative. Hence
\[K_{\ga}(X,\xi)\,K_{\ga}(X,Z)\geq R_{\ga}(X,Z,\xi,X)^2.\]
If $(M,\ga)$ has totally geodesic fibers, Gray--O'Neill's equation \cite{oneill,gw} gives $R_{\ga}(X,Z,\xi,X)=\lr{(\nabla_XA^*)_X\xi,Z}$ and $K_{\ga}(\xi,X)=|A^*_X\xi|^2$. On a relatively compact neighborhood there is $C>0$ with $K_{\ga}(X,Z)\leq C|X|^2|Z|^2$, whence
\[\lr{(\nabla_XA^*)_X\xi,Z}=R_{\ga}(X,Z,\xi,X)\leq \sqrt{K_{\ga}(X,\xi)\,K_{\ga}(X,Z)}\leq \sqrt C\,|X||Z|\,|A^*_X\xi|,\]
which is \eqref{eq:WNN-intro} with $\tau=\sqrt C$. \qedhere
\end{proof}

\begin{proof}[Proof of item \ref{item:essevainofim}]
We argue by contradiction. Suppose there are non-zero $X\in\cal H$, $\nu_0\in\cal V$ with $A^*_X\nu_0=0$. Lemma \ref{lem:flatgeo} gives $A^*_{\dot c}\nu(t)=0$ for the dual holonomy field $\nu$ along $c(t)=\exp(tX)$ with $\nu(0)=\nu_0$. Since $\sec(M)>0$, $\pi$ has positive vertizontal curvature, thus, it is integrally fat. However,  $|A^*_{\dot c}u(t)|\equiv 0$ along $c$ for $u=\nu/|\nu|$, contradicting Definition~\ref{dfn:integral-fatness}. \qedhere
\end{proof}

\vspace{1em}

\section{\texorpdfstring{$A$-Clifford algebras on Riemannian submersions}{A-Clifford algebras on Riemannian submersions}}
\label{sec:aclifford}

Recall the following definition, first appearing in \cite{MOROIANU2011940}:
\begin{definition}\label{def:clifford-structure}
A \emph{rank \texorpdfstring{r}{r} Clifford structure} on a Riemannian manifold \texorpdfstring{$(B, b)$}{(B, b)} is an oriented
rank \( r \) Euclidean bundle \( (E,\ga_E) \) over \( B \), together with a non-vanishing algebra bundle morphism,
called the \emph{Clifford morphism},
\[
\varphi : \mathrm{Cl}(E,\ga_E) \to \mathrm{End}(TB),
\]
which maps \( E \) into the bundle of skew-symmetric endomorphisms \( \mathrm{End}^{-}(TB) \).
\end{definition}

In this section, we aim at combining Clifford structures with Riemannian submersions and explore the curvature properties of the resulting structure.

\begin{definition}\label{def:a-clifford}
A Riemannian submersion $\pi:(M,\ga)\to(B,b)$ is \emph{$A$-Clifford} if $(A^u)^2$ is a multiple of the identity for every $u\in\cal V$. We then write
\begin{equation}\label{eq:Ndef}
(A^u)^2=-N(u)\,\mathrm{Id}_{\cal H},\qquad N(u)\ge 0 .
\end{equation}
It is \emph{$A$-Clifford in the normalized gauge} if, moreover, $N(u)=|u|^2_{\ga}$; i.e., \ if $u\mapsto A^u$ is a Clifford morphism for $\mathrm{Cl}(\cal V,\ga|_{\cal V})$ in the sense of Definition \ref{def:clifford-structure}.
\end{definition}

For short, we shall call $A$-Clifford submersions in the normalized gauge as normalized $A$-Clifford submersions. 

As we shall see in this section, (normalized) $A$-Clifford Riemannian submersions have constrained curvature identities. In particular, they are the \emph{extreme case of weak nonnegative curvature}: the right-hand side of the WNN inequality \eqref{eq:WNN-intro} vanishes under normalized $A$-Clifford submersions with totally geodesic fibers.

\begin{theorem}\label{thm:aclifford-basics}
Let $\ga$ be adapted to $\pi:(M,\cal H)\to(B,b)$ and $\dim\cal V=r$. Then $\ga$ is $A$-Clifford if and only if every metric adapted to $(\cal H, b)$ is $A$-Clifford. Moreover:
\begin{enumerate}
\item\label{item:acl-N} $N$ in \eqref{eq:Ndef} is a quadratic form on $\cal V$, and
\[|A^*_Xu|^2=N(u)|X|^2\qquad\text{for all }X\in\cal H,\ u\in\cal V.\]
In particular, $\pi$ is fat if and only if $N$ is definite;
\item\label{item:acl-gauge} if $\pi$ is fat, then $N$ is the unique vertical metric, among those adapted to $(\cal H, b)$, for which the structure is in the normalized gauge. In this case $u\mapsto A^u$ is a rank $r$ Clifford structure over $(B,b)$ and $|A^*_Xu|^2=|u|^2|X|^2$; if, in addition, the fibers are totally geodesic in the normalized gauge, then
\[K(X,u)=|A^*_Xu|^2=|u|^2|X|^2 ;\]
\item\label{item:acl-wnn} if $\pi$ is fat, then $\ga$ is weakly nonnegatively curved;
\item\label{item:acl-tau} if $\pi$ is $A$-Clifford in the normalized gauge, then the WNN condition holds taking $\tau=\sigma$, the maximum operator norm of $S_X$ over unit $X\in\cal H$. Moreover, the right-hand side of \eqref{eq:WNN-intro} vanishes identically when the fibers are totally geodesic, and $N(\nu)=|\nu|^2$ is constant on holonomy fields.
\end{enumerate}
\end{theorem}

\begin{proof}
Let $\ga'$ be adapted to $(\cal H, b)$ and let $P$ be as in Lemma \ref{lem:dualinv}, so that $A^{\dagger}_Xu=A^*_X Pu$. Then, the endomorphism attached to $u$ by $\ga'$ is $(A^{\dagger})^u=A^{Pu}$. As $P$ is an isomorphism of $\cal V_p$, the families $\{A^u\}_{u\in\cal V_p}$ and $\{(A^{\dagger})^u\}_{u\in\cal V_p}$ coincide as subsets of $\End(\cal H_p)$. Since all adapted metrics agree on $\cal H$, being a multiple of the identity is a metric-independent property of these endomorphisms.

Item \ref{item:acl-N} follows by polarizing $(A^{u+v})^2=-N(u+v)\mathrm{Id}$ and using $(A^u)^2=-N(u)\mathrm{Id}$, $(A^v)^2=-N(v)\mathrm{Id}$ gives
\[A^uA^v+A^vA^u=-\bigl(N(u+v)-N(u)-N(v)\bigr)\mathrm{Id},\]
whose left-hand side is bilinear in $(u,v)$. Hence, $N$ is the quadratic form defined by the symmetric bilienear form $N(u,v) = \frac{1}{2}[N(u+v)-N(u)-N(v)]$. Moreover,
\[|A^*_Xu|^2=\ga(A^uX,A^uX)=-\ga((A^u)^2X,X)=N(u)|X|^2 .\]
Fatness means $A^*_Xu\ne0$ for all $X\ne0\ne u$, which is translated to $N(u)>0$ for $u\neq 0$.

For item \ref{item:acl-gauge}, note that $N$ defines a  vertical metric $|\cdot|_N$ whenever $\pi$ is fat. By construction, $|(A^u)^2X|^2=|-N(u)X|^2=|u|^2_N|X|^2$, which is the normalization of Definition \ref{def:clifford-structure} for $\mathrm{Cl}(\cal V,N)$. Thus, $u\mapsto A^u$ extends to an algebra bundle morphism. Uniqueness is straightforward. As for the curvature, Gray--O'Neill gives $K(X,u)=|A^*_Xu|^2-|S_Xu|^2+\ga((\nabla^{v}_X S)_Xu,u)$, which reduces to $K(X,u)=|A^*_Xu|^2$ when $S\equiv0$; combined with $|A^*_Xu|^2=N(u)|X|^2=|u|^2|X|^2$ this is the displayed identity. 

For proving item \ref{item:acl-wnn}, recall that the right-hand side of \eqref{eq:WNN-intro} equals $\tfrac12 \frac{\mathrm d}{\mathrm dt}|A^*_{\dot c}\nu|^2$ along the unit speed geodesic generated by $X$, where $\nu$ is the dual holonomy field with $\nu(0)=\xi$. By \ref{item:acl-N}, $|A^*_{\dot c}\nu(t)|^2=N(\nu(t))$, so, denoting by $N(\cdot,\cdot)$ the polarization and using $\nabla^{v}_{\dot c}\nu=S_{\dot c}\nu$,
\begin{equation}\label{eq:Nderivative}
\frac{\mathrm d}{\mathrm dt}N(\nu,\nu)=(\nabla_{\dot c}N)(\nu,\nu)+2N(S_{\dot c}\nu,\nu)\le \bigl(|\nabla N|+2\sigma|N|\bigr)|\nu|^2 .
\end{equation}
Here, $\nabla N$ denotes the derivative of the $(0,2)$-tensor $N$ on $\cal V$ with respect to the projected connection $D_X:=(\nabla_X\cdot)^{v}$, and we use $D_{\dot c}\nu=\nabla^{v}_{\dot c}\nu=S_{\dot c}\nu$ from \eqref{eq:dualhol}. In a relatively compact neighborhood, the tensor norms $|\nabla N|$ and $|N|$ are bounded, and $N$ being positive definite by fatness, there is $c>0$ with $N(\nu,\nu)\ge c|\nu|^2$. Hence, the right-hand side of \eqref{eq:Nderivative} is at most $2\tau N(\nu,\nu)$ with $2\tau=(|\nabla N|+2\sigma|N|)/c$. Since $N(\nu,\nu)=|A^*_{\dot c}\nu|^2$, we verify \eqref{eq:WNN-intro}.

\ref{item:acl-tau} Notice that the projected connection $D_X=(\nabla_X\cdot)^{v}$ is metric for $\ga|_{\cal V}$. Hence, $N=\ga|_{\cal V}$ satisfies $\nabla N=0$. Thus, \eqref{eq:Nderivative} reduces to
\[\frac{\mathrm d}{\mathrm dt}N(\nu,\nu)=\frac{\mathrm d}{\mathrm dt}|\nu|^2=2\ga(S_{\dot c}\nu,\nu)\le2\sigma|\nu|^2.\]
The above equation combined with Equation \eqref{eq:mixed-curvature-is-a-norm} implies that $\ga$ satisfies \eqref{eq:WNN-intro} with $\tau=\sigma$. If, moreover, $S\equiv0$, the same identity gives $\frac{\mathrm d}{\mathrm dt}|\nu|^2=0$. Thus, the right-hand side of \eqref{eq:WNN-intro} vanishes identically and $N(\nu)=|A^*_{\dot c}\nu|^2=|\nu|^2$ is constant along the holonomy field $\nu$. \qedhere
\end{proof}

\begin{theorem}\label{thm:adapted-totally-geodesic}
    Let $\pi:(M,\ga)\rightarrow (B,b)$ be a Riemannian submersion with totally geodesic fibers. Then, $\pi$ is $A$-Clifford in the normalized gauge if, and only if, for all horizontal $X,Y$ and all vertical $\xi,\eta$,
    \begin{equation}\label{eq:curvature-characterization}
        R_\ga(\xi,Y,X,\eta)+R_\ga(\eta,Y,X,\xi)=2\,\ga(X,Y)\,\ga(\xi,\eta).
    \end{equation}
    If, in addition, $R_\ga(\xi,Y,X,\eta)$ is symmetric in $(\xi,\eta)$, then \eqref{eq:curvature-characterization} is equivalent to
    \begin{equation}\label{eq:analogously}
        R^{\mathbf{v}}(\xi,Y)X=\ga(X,Y)\xi .
    \end{equation}
   In either case, $\ga$ has vertizontal curvature constant equal to $1$ on non-degenerate planes $X\wedge\xi$.
\end{theorem}
\begin{proof}
First, notice that for any Riemannian submersion $\pi$ with totally geodesic fibers, an easy adaptation of Lemma 2.17 in \cite{Ziller_fatnessrevisited} gives for basic vector fields $X, Z$ and holonomy fields $\xi_a,\xi_b$ along the geodesic generated by $X$,
    \begin{equation}\label{eq:ziller217}
   \ga(J_a\circ J_b(X),Z)=\ga(J_{[a,b]}(X),Z)-R_{\ga}(\xi_a,Z,X,\xi_b),
    \end{equation}
  where $J_a:=-A^{\xi_a}$, $J_b:=-A^{\xi_b}$ and $J_{[a,b]}:=-A^{[\xi_a,\xi_b]}$. We shall only use the following consequence of \eqref{eq:ziller217}: symmetrizing in $a,b$ and using $[\xi_a,\xi_b]=-[\xi_b,\xi_a]$, so that the $J_{[a,b]}$-contributions cancel,
\[\ga\bigl((J_aJ_b+J_bJ_a)(X),Z\bigr)=-R_\ga(\xi_a,Z,X,\xi_b)-R_\ga(\xi_b,Z,X,\xi_a).\]
Polarizing \eqref{eq:Ndef} with $N=\ga|_{\cal V}$ we check that the normalized gauge $A$-Clifford condition is $J_aJ_b+J_bJ_a=-2\ga(\xi_a,\xi_b)\mathrm{Id}$, which is exactly \eqref{eq:curvature-characterization}. Under the additional symmetry, \eqref{eq:curvature-characterization} halves to $R_\ga(\xi,Y,X,\eta)=\ga(X,Y)\ga(\xi,\eta)$ for all vertical $\eta$. Since $\ga|_{\cal V}$ is nondegenerate this is \eqref{eq:analogously}. Finally, with totally geodesic fibers, $K(X,\xi)=|A^*_X\xi|^2=|\xi|^2|X|^2$ by Theorem \ref{thm:aclifford-basics}\ref{item:acl-gauge}.
\end{proof}
Recall from \cite{boyer1993sasakian} that a Riemannian manifold \( (P^m, \ga) \) is \emph{$3$-Sasakian} if it carries three Killing fields $\xi_1,\xi_2,\xi_3$ with
\begin{equation}\label{eq:3sas}
\ga(\xi_a,\xi_b)=\delta_{ab},\qquad [\xi_a,\xi_b]=2\varepsilon_{abc}\xi_c,\qquad (\nabla_X\phi_a)Y=\ga(X,Y)\xi_a-\eta_a(Y)X,
\end{equation}
where $\phi_a=-\nabla\xi_a$ and $\eta_a=\ga(\xi_a,\cdot)$. Locally the $\xi_a$ generate an isometric $G=\mathrm{SO}(3)$ (or $\mathrm S^3$) action. When this action is free, $P$ is the total space of a $G$-principal bundle, which we refer to as a \emph{$3$-Sasakian principal bundle}.

We record the two consequences of \eqref{eq:3sas} used below. Differentiating $\nabla_Y\xi_a=-\phi_aY$ and using the third identity in \eqref{eq:3sas},
\[
R(X,Y)\xi_a=-(\nabla_X\phi_a)Y+(\nabla_Y\phi_a)X=\eta_a(Y)X-\eta_a(X)Y ,
\]
so that, for horizontal $X,Y$ and vertical $V$,
\begin{equation}\label{eq:3sas-mixed}
R(X,Y)V=0 ,
\end{equation}
and, taking $Y=\xi_b$ and pairing, for horizontal $X,Y$ and vertical $\xi,\eta$,
\begin{equation}\label{eq:3sas-vhhv}
R_\ga(\xi,Y,X,\eta)=\ga(X,Y)\,\ga(\xi,\eta).
\end{equation}
Moreover $\nabla_{\xi_a}\xi_b=\h[\xi_a,\xi_b]=\varepsilon_{abc}\xi_c$ is vertical, so the orbits are totally geodesic, and they carry a bi-invariant metric of curvature $K(\xi_1,\xi_2)=\tfrac14|[\xi_1,\xi_2]|^2=1$.

\begin{corollary}\label{cor:A-clifford-Sasakian}
    Every $3$-Sasakian principal bundle $(P,\ga)\rightarrow (B,b)$ has totally geodesic fibers of curvature $1$ and is $A$-Clifford in the normalized gauge.
\end{corollary}
\begin{proof}
The fibers are the orbits, which are totally geodesic of curvature $1$ by the discussion above. Identity \eqref{eq:3sas-vhhv} is \eqref{eq:analogously}; symmetrizing it in $(\xi,\eta)$ yields \eqref{eq:curvature-characterization}, and Theorem \ref{thm:adapted-totally-geodesic} gives that $\pi$ is $A$-Clifford in the normalized gauge.
\end{proof}

Next, we take a new look at a theorem of Dearricott. We recall a lemma from \cite{CDR} and show how this simplifies the curvature of a 3-Sasakian manifold.

\begin{theorem}[Theorem 1 in \cite{Dearricott2004}]\label{cor:dearrictott}
    A $3$-Sasakian principal bundle $\pi: (P,\ga)\rightarrow (B,b)$ admits a metric with positive sectional curvature if, and only if, $(B,b)$ has positive sectional curvature.
\end{theorem}

That follows, essentially, from the normal form of \cite[Lemma 1]{CDR}; we include the short argument.

\begin{lemma}\label{lem:cdr-normal-form}
Let $\pi:(M,\ga)\to(B,b)$ be a Riemannian submersion and $\sigma\subset T_pM$ a plane. Then $\sigma$ admits a $\ga$-orthonormal basis $\{X_1+V_1,\ X_2+V_2\}$, with $X_i\in\cal H_p$ and $V_i\in\cal V_p$, such that
\[\ga(X_1,X_2)=0\qquad\text{and}\qquad \ga(V_1,V_2)=0 .\]
\end{lemma}
\begin{proof}
The assignment $q(v,w):=\ga(v^{h},w^{h})$ is a symmetric bilinear form on the two-dimensional space $\sigma$, so it is diagonalized by some $\ga|_\sigma$-orthonormal basis $\{e_1,e_2\}$ of $\sigma$. Writing $e_i=X_i+V_i$ with $X_i=e_i^{h}$ and $V_i=e_i^{v}$, diagonality reads $\ga(X_1,X_2)=q(e_1,e_2)=0$; and since $0=\ga(e_1,e_2)=\ga(X_1,X_2)+\ga(V_1,V_2)$, also $\ga(V_1,V_2)=0$.
\end{proof}

\begin{proposition}\label{prop:3sas-decoupling}
Let $\pi:(P,\ga)\to(B,b)$ be a $3$-Sasakian principal bundle, $\sigma\subset T_pP$ a plane, and $\{X_1+V_1,X_2+V_2\}$ a normal-form basis of $\sigma$ as in Lemma \ref{lem:cdr-normal-form}. Then the sectional curvature of $\sigma$ decouples as
\begin{equation}\label{eq:decoupling}
\mathrm{sec}_\ga(\sigma)=K_\ga(X_1,X_2)+|A^*_{X_2}V_1|^2+|A^*_{X_1}V_2|^2+|V_1|^2|V_2|^2 ,
\end{equation}
where $K_\ga(X_1,X_2)=K_b(X_1,X_2)-3|A_{X_1}X_2|^2$ by Gray--O'Neill's horizontal curvature equation.
\end{proposition}
\begin{proof}
Since the given basis is orthonormal, $\mathrm{sec}_\ga(\sigma)=R_\ga(X_1+V_1,X_2+V_2,X_2+V_2,X_1+V_1)$. Expanding by multilinearity produces sixteen terms, which we group according to the number of vertical entries.

\emph{No vertical entry.} The single term is $R_\ga(X_1,X_2,X_2,X_1)=K_\ga(X_1,X_2)$.

\emph{One vertical entry.} There are four such terms, and each is, after using the antisymmetry in the first two and in the last two arguments together with the pair symmetry, of the form $\pm\ga(R(X,Y)V,Z)$ with $X,Y,Z$ horizontal and $V$ vertical. All vanish by \eqref{eq:3sas-mixed}.

\emph{Two vertical entries.} There are six such terms. Two of them, $R_\ga(V_1,X_2,X_2,V_1)$ and $R_\ga(X_1,V_2,V_2,X_1)$, are the unreduced vertizontal curvatures $K_\ga(V_1,X_2)$ and $K_\ga(X_1,V_2)$; since the fibers are totally geodesic these equal $|A^*_{X_2}V_1|^2$ and $|A^*_{X_1}V_2|^2$. Two further terms, $R_\ga(V_1,V_2,X_2,X_1)$ and $R_\ga(X_1,X_2,V_2,V_1)$, equal $\pm\ga(R(X_1,X_2)V_i,V_j)$ after the pair symmetry, hence vanish by \eqref{eq:3sas-mixed}. The remaining two are
\[
R_\ga(V_1,X_2,V_2,X_1)=-R_\ga(V_1,X_2,X_1,V_2),\qquad
R_\ga(X_1,V_2,X_2,V_1)=-R_\ga(V_2,X_1,X_2,V_1),
\]
and by \eqref{eq:3sas-vhhv} both equal $-\ga(X_1,X_2)\ga(V_1,V_2)$, which vanishes because the basis is in normal form.

\emph{Three vertical entries.} There are four such terms, each of the form $R_\ga(U,V,W,X)$ with $U,V,W$ vertical and $X$ horizontal, up to sign. By Gray--O'Neill's equations \cite[p.44]{gw} these are expressed solely through the second fundamental form of the fibers, hence vanish.

\emph{Four vertical entries.} The single term $R_\ga(V_1,V_2,V_2,V_1)$ equals the corresponding curvature of the fiber, the fibers being totally geodesic; as the fibers have curvature $1$ and $\ga(V_1,V_2)=0$, it equals $|V_1|^2|V_2|^2$.

Summing the surviving terms gives \eqref{eq:decoupling}.
\end{proof}

Recall that $\mathrm{Cl}_{0,3}\cong\bb H\oplus\bb H$ (here $\bb H$ stands for the quaternions) has exactly two inequivalent irreducible modules, distinguished by the action of the volume element $\omega=e_1e_2e_3$, which is central and squares to $\mathrm{Id}$: it acts as $+\mathrm{Id}$ on one and as $-\mathrm{Id}$ on the other. A $\mathrm{Cl}_{0,3}$-module is \emph{isotypic} if only one of them occurs. Equivalently, if $A^{\xi_1}A^{\xi_2}=\pm A^{\xi_3}$ with a constant sign. Lemma \ref{lem:3sas-isotypic} below shows that for $3$-Sasakian principal bundles the module is always isotypic, with $\omega=-\mathrm{Id}$. Theorem \ref{thm:3sasakian-characterization} provides the converse. Together with Corollary \ref{cor:A-clifford-Sasakian} yields Theorem \ref{ithm:3sas}. The conclusion that $(B,b)$ is quaternionic K\"ahler with positive scalar curvature is \cite[Proposition 2.16]{Ziller_fatnessrevisited}.

\begin{lemma}\label{lem:3sas-isotypic}
Let $\pi:G\hookrightarrow (P,\ga)\to(B,b)$ be a $3$-Sasakian principal bundle with $G=\mathrm S^3$ or $\mathrm{SO}(3)$. Then, for every $p\in P$, the endomorphisms $A^{\xi_a}$ of $\cal H_p$ satisfy
\begin{equation}\label{eq:quaternionic-relation}
A^{\xi_a}A^{\xi_b}=\varepsilon_{abc}A^{\xi_c}-\delta_{ab}\,\mathrm{Id}_{\cal H_p} .
\end{equation}
In particular $A^{\xi_1}A^{\xi_2}=A^{\xi_3}$, the volume element acts by $\omega:=A^{\xi_1}A^{\xi_2}A^{\xi_3}=-\mathrm{Id}$, and the $\mathrm{Cl}_{0,3}$-module $\cal H_p$ is isotypic.
\end{lemma}
\begin{proof}
Write $\phi_a=-\nabla\xi_a$ as in \eqref{eq:3sas}, so that $\phi_aW=-\nabla_W\xi_a$ for every $W$.
 For horizontal $X$, the Killing equation for $\xi_a$ gives
\[\ga(\nabla_X\xi_a,\xi_b)=-\ga(\nabla_{\xi_b}\xi_a,X)=0.\]
 Thus, $\nabla_X\xi_a$ is horizontal, so $\phi_a$ maps $\cal H$ into $\cal H$. Hence,
\begin{equation}\label{eq:vai-igualar}A^{\xi_a}X=A^*_X\xi_a=-(\nabla_X\xi_a)^{h}=-\nabla_X\xi_a=\phi_aX,\qquad X\in\cal H .\end{equation}

Evaluation $Y=\xi_b$ and horizontal $X$ on \eqref{eq:3sas}, its right-hand side is
\[\ga(X,\xi_b)\xi_a-\eta_a(\xi_b)X=-\delta_{ab}X.\]
 Its left-hand side expands as
\[(\nabla_X\phi_a)\xi_b=\nabla_X(\phi_a\xi_b)-\phi_a(\nabla_X\xi_b)
=\varepsilon_{abc}\nabla_X\xi_c-\phi_a(-\phi_bX)
=-\varepsilon_{abc}\phi_cX+\phi_a\phi_bX .\]
Equating the two expressions gives $\phi_a\phi_b=\varepsilon_{abc}\phi_c-\delta_{ab}\mathrm{Id}$ on $\cal H$, which is \eqref{eq:quaternionic-relation} by \eqref{eq:vai-igualar}.
Taking $a=b$ in \eqref{eq:quaternionic-relation} yields $(A^{\xi_a})^2=-\mathrm{Id}$, and taking $(a,b,c)=(1,2,3)$ yields $A^{\xi_1}A^{\xi_2}=A^{\xi_3}$. Therefore
\[\omega=A^{\xi_1}A^{\xi_2}A^{\xi_3}=(A^{\xi_3})^2=-\mathrm{Id}.\]
\end{proof}

\begin{theorem}\label{thm:3sasakian-characterization}
Let $\pi: G\hookrightarrow (P,\ga)\to(B,b)$ be a Riemannian principal bundle with $G=\mathrm S^3$ or $\mathrm{SO}(3)$ and totally geodesic fibers of curvature $1$, and assume that $\pi$ is $A$-Clifford in the normalized gauge. Let $\xi_1,\xi_2,\xi_3$ be an orthonormal basis of Killing fields generating the $G$-action, normalized so that $[\xi_a,\xi_b]=2\varepsilon_{abc}\xi_c$, and let
\[\omega:=A^{\xi_1}A^{\xi_2}A^{\xi_3}\]
be the associated volume element. If $\omega=-\mathrm{Id}$ (equivalently, if $A^{\xi_1}A^{\xi_2}=A^{\xi_3}$) then $(P,\ga)$ is $3$-Sasakian.
\end{theorem}
\begin{proof}
Set $\phi_a:=-\nabla\xi_a$ and $\eta_a:=\ga(\xi_a,\cdot)$. We check that these data assemble into a contact $3$-structure on $P$. By Kashiwada's theorem \cite{Kashiwada2001} it follows that it defines a Sasakian $3$-structure on $P$. 

As for the first step, we check that each $(\phi_a,\eta_a)$ defines a contact structure on $P$: Since the fibers are totally geodesic, then $\phi_a|_{\cal H}=A^{\xi_a}$. On $\cal H$, the normalized gauge gives $\phi_a^2=(A^{\xi_a})^2=-\mathrm{Id}$, and on $\cal V$ one computes $\phi_a^2\xi_a=0$ and $\phi_a^2\xi_b=-\xi_b$ for $b\ne a$. Hence, $\phi_a^2=-\mathrm{Id}+\eta_a\otimes\xi_a$. Moreover, $\ga(\phi_a\cdot,\phi_a\cdot)=\ga-\eta_a\otimes\eta_a$ since $\ga(A^{\xi_a}X,A^{\xi_a}Y)=-\ga((A^{\xi_a})^2X,Y)=\ga(X,Y)$. For horizontal $X,Y$,
\begin{align}
    \mathrm d\eta_a(X,Y)&=-\eta_a([X,Y])=-2\ga(A_XY,\xi_a)=-2\ga(A^{\xi_a}X,Y)=2\ga(X,\phi_aY)\\
    \mathrm d\eta_a(\xi_b,\xi_c)&=-\eta_a([\xi_b,\xi_c])=-2\varepsilon_{bca}=2\ga(\xi_b,\phi_a\xi_c),
\end{align}
as required \cite[p. 830]{Kashiwada2001}.

The proof is concluded if one verifies
\[\phi_c=\phi_a\phi_b-\eta_b\otimes\xi_a\qquad\text{for cyclic }(a,b,c).\footnote{With the conventions of \cite{Kashiwada2001} the relation is written $\phi_c=\phi_a\phi_b-\eta_a\otimes\xi_b$; the form used here is the same relation for our normalization $[\xi_a,\xi_b]=2\varepsilon_{abc}\xi_c$ and $\phi_a=-\nabla\xi_a$.}\]
On $\cal V$, taking $(a,b,c)=(1,2,3)$: at $\xi_1$ one has $\phi_1\phi_2\xi_1=\phi_1(-\xi_3)=\xi_2$ and $\eta_2(\xi_1)\xi_1=0$, matching $\phi_3\xi_1=\xi_2$; at $\xi_2$, $\phi_1\phi_2\xi_2=0$ and $\eta_2(\xi_2)\xi_1=\xi_1$, matching $\phi_3\xi_2=-\xi_1$; at $\xi_3$, $\phi_1\phi_2\xi_3=\phi_1\xi_1=0$ and $\eta_2(\xi_3)=0$, matching $\phi_3\xi_3=0$. The other cyclic permutations are checked in the same way. On $\cal H$ the $\eta$-term vanishes and the identity reads $A^{\xi_a}A^{\xi_b}=A^{\xi_c}$, which is exactly the hypothesis: the Clifford relations of Theorem \ref{thm:aclifford-basics} give $A^{\xi_a}A^{\xi_b}=-A^{\xi_b}A^{\xi_a}$ for $a\neq b$ and $(A^{\xi_a})^2=-\mathrm{Id}$, so $\omega$ is central with $\omega^2=\mathrm{Id}$ and $(A^{\xi_3})^{-1}=-A^{\xi_3}$. Hence, $A^{\xi_1}A^{\xi_2}=\omega\,(A^{\xi_3})^{-1}=-\omega A^{\xi_3}$, and $\omega=-\mathrm{Id}$ if and only if $A^{\xi_1}A^{\xi_2}=A^{\xi_3}$. The cyclic permutations follows by the same computation.
\end{proof}


\vspace{1em}

\section*{Acknowledgments}
Some ideas explored in this work were developed during the first-named author's postdoctoral position at the University of Fribourg, supported partly by the SNSF Project 200020E\_193062 and the DFG Priority Program SPP 2026. Part of this work was also conceived when L. F. Cavenaghi received support from CAPES and The São Paulo Research Foundation (FAPESP), with grants 2022/09603-9 and 2023/14316-1. 

L. Grama is partially supported by FAPESP grants 2018/13481-0, 2021/04065-6, 2023/13131-8, and CNPq grant 306021/2024-2.

L F. Cavenaghi is supported by the Simons Foundation, grant SFI-MPS-T-Institutes-00007697, and the Ministry of Education and Science of the Republic of Bulgaria, grant DO1-239/10.12.2024.

\bibliographystyle{plain}
	\bibliography{main}

	\end{document}